\documentclass[12pt]{amsart}
\usepackage{enumerate}
\usepackage{enumitem}
\usepackage[colorlinks=true, linkcolor=blue, urlcolor=blue, citecolor=blue, anchorcolor=blue, pdfborder={0 0 0}]{hyperref}
\usepackage{url}
\usepackage{graphicx,color,colortbl}
\usepackage[dvipsnames,table,xcdraw]{xcolor}
\usepackage{cite}
\usepackage{amsthm, amsmath, amssymb}
\usepackage{mathtools}
\usepackage[top=45truemm, bottom=45truemm, left=30truemm, right=30truemm]{geometry}
\usepackage{subcaption}
\usepackage{cancel}
\usepackage{float}
\usepackage{tabularx}
\usepackage{makecell}
\usepackage{array}
\usepackage{ragged2e}
\usepackage{csquotes}
\usepackage{booktabs}

\usepackage{array}
\usepackage{threeparttable} 

\definecolor{LightGray}{RGB}{242,242,242}

\usepackage{listings}
\newcolumntype{P}[1]{>{\RaggedRight\hspace{0pt}}p{#1}}
\newcolumntype{L}{>{\begin{math}}l<{\end{math}}}%
\newcolumntype{C}{>{\begin{math}}c<{\end{math}}}%
\newcolumntype{R}{>{\begin{math}}r<{\end{math}}}%

\title[Hearing the Sides: Recovering a Planar Rectangle from Eigenvalues]{Hearing the Sides: Recovering a Planar Rectangle from Eigenvalues}

\author[Eldar Sultanow]{\href{https://orcid.org/0000-0001-5257-2236}{\includegraphics[scale=0.06]{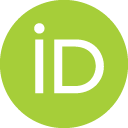}}\hspace{1mm}Eldar Sultanow}
\address{Eldar Sultanow\\Capgemini Deutschland GmbH\\Nuremberg, Germany}
\curraddr{}
\email{eldar.sultanow@capgemini.com}

\author[A.\ Hatziiliou]{Andreas Hatziiliou}
\address{Andreas Hatziiliou\\University of British Columbia, Department of Mathematics\\1984 Mathematics Rd, Vancouver BC\\Canada}
\curraddr{}
\email{ahatzi@math.ubc.ca}


\begin{document}

\begingroup
\let\MakeUppercase\relax
\maketitle
\endgroup

\begin{abstract}
We present a direct, index-free method to recover the side lengths of a planar rectangle the spectrum of its Dirichelet Laplacian, assuming only access to a finite subset of eigenvalues. No modal indices $(m,n)$ are available, and the list may begin at an arbitrary unknown offset; in particular, the lowest eigenvalues may be missing, so classical formulas based on $\lambda_{1,0}$ and $\lambda_{0,1}$ cannot be used. Our reconstruction procedure extracts geometric information solely from the asymptotic density and oscillatory structure of the ordered spectrum. The area $ab$ is obtained from the high-frequency Weyl slope, while the fundamental lengths $2a$ and $2b$ appear as dominant periodic--orbit contributions in the Fourier transform of the spectral fluctuations. This separation of smooth and oscillatory components yields a robust, offset-agnostic recovery of both side lengths. The result is a fully index-free algorithm that reconstructs the geometry of a rectangular planar domain even when the spectrum is incomplete and all modal information is lost.
\end{abstract}

\section{The Problem}
Inverse spectral problems of this form have been studied for decades, from Kac’s classical question to modern analyses of simple planar domains such as the triangle by Chang and DeTurck \cite{ChangDeTurck1989}. Recovering geometric data from the Dirichlet spectrum remains an active topic in spectral geometry and mathematical physics.

The Laplace spectrum of a planar domain carries rich geometric information. For highly symmetric regions such as discs or rectangles, the corresponding eigenvalue problems can be solved explicitly, and the structure of the spectrum reflects the underlying geometry in a remarkably direct way. This raises a natural inverse question: To what extent can the geometry of such a domain be recovered from its spectrum alone? In this paper we focus on the simplest nontrivial case: a rectangular domain \cite[p.~161]{Strauss2008}
\[
\Omega_{a,b} = (0,a)\times(0,b)\subset\mathbb{R}^2,
\qquad a,b>0,
\]
equipped with Dirichlet boundary conditions. Separation of variables yields an explicit formula for the eigenvalues of the Dirichlet Laplacian, see \cite[p.~322]{Strauss2008}:
\begin{equation}
\label{eq:rect-spectrum}
\lambda_{m,n}(a,b)=\left(\frac{m\pi}{a}\right)^2+\left(\frac{n\pi}{b}\right)^2,\qquad m,n\in\mathbb{N},
\end{equation}
so the full spectrum consists of all pairwise sums of the one-dimensional Dirichlet eigenvalues on intervals of lengths $a$ and $b$. Equation~\eqref{eq:rect-spectrum} completely characterises the spectrum: $\sigma(-\Delta_{\Omega_{a,b}})$ is the multiset of all values on the right-hand side.

The direct problem is therefore elementary. The inverse problem, however, is considerably more subtle. Throughout we assume that only the \emph{ordered} sequence of Dirichlet eigenvalues is known, written in nondecreasing order as in \cite[p.~310]{Strauss2008}, 
\cite[p.~244]{PinchoverRubinstein2005}:
\begin{equation}
\label{eq:enumerated-eigenvalues}
0 < \lambda_1 \le \lambda_2 \le \lambda_3 \le \cdots,
\end{equation}
With no access to the underlying mode indices $(m,n)$, can one reconstruct the side lengths $a$ and $b$?

\section{The Idea}
The central observation is that the geometric data of the rectangle can be recovered directly from the ordered Dirichlet spectrum, without access to the underlying indices $(m,n)$ arising from the separation of variables on the rectangle, which label the eigenmodes $\sin(m\pi x/a)\cdot\sin(n\pi y/b)$, see \cite[p.~322]{Strauss2008}, and generate the eigenvalues $\lambda_{m,n}$ in \eqref{eq:rect-spectrum}. In the inverse problem, however, these labels are not available; only the nondecreasing enumeration~\eqref{eq:enumerated-eigenvalues} is given.

Our approach is to extract the geometric invariants $ab$, $a$, and $b$ directly from this ordered sequence. Define the eigenvalue counting function, see e.g. \cite[p.~17]{Levitin_etal2023}, by
\begin{equation}
\label{eq:eigenvalue-counting-function}
N(\lambda)=\#\{k : \lambda_k \le \lambda\}.
\end{equation}
The counting function admits a decomposition into a smooth Weyl term and a remainder term (fluctuations):
\begin{equation}
\label{eq:decomposition-smooth-osc}
N(\lambda)=N_{\mathrm{smooth}}(\lambda)+F(\lambda),
\end{equation}
which is completely analogous to the asymptotic formula given, for example, by Levitin et al.~\cite[p.~18]{Levitin_etal2023}:
\[
N(\lambda)=\frac{\mathrm{Area}(T_a^2)}{4\pi}\lambda + R(\lambda),
\]
where $R(\lambda)=o(\lambda)$ as $\lambda\to\infty$, in the sense described by de~Bruijn, see Appendix~\ref{appx:weyls-law}. This decomposition allows us to isolate the smooth part governing the area $ab$ and the oscillatory remainder encoding the periodic--orbit lengths, thereby providing a complete reconstruction of $a$ and $b$ from the ordered Dirichlet spectrum alone. The geometric-spectral duality underlying this split (formalized through Poisson summation and, in the integrable case, the Berry-Tabor theory, which shows that the oscillatory part of the level density of an integrable system can be expressed in terms of its classical periodic orbits\cite[p.~89]{BrackBhaduri1997}) relates the oscillatory remainder of the counting function directly to the lengths of periodic classical orbits on the domain.

\section{Recovering the Area}
The first geometric invariant encoded in the ordered Dirichlet spectrum of $\Omega_{a,b}$ is the \emph{area} $ab$. This information is contained in the high-frequency growth of the eigenvalue counting function~\eqref{eq:eigenvalue-counting-function}, which is accessible from the enumerated spectral sequence $\{\lambda_k\}_{k\in\mathbb{N}}$. For planar domains, Weyl's law, as recalled in Appendix~\ref{appx:weyls-law}, asserts that
\begin{equation}
\label{eq:rectangle-weyl}
N(\lambda)
= \frac{\operatorname{Area}(\Omega_{a,b})}{4\pi}\,\lambda + o(\lambda)
= \frac{ab}{4\pi}\,\lambda + o(\lambda),
\qquad \lambda\to\infty.
\end{equation}
Here the symbol $o(\lambda)$ is the standard \emph{little-$o$} notation from asymptotic analysis (not an error function): it denotes any function $f(\lambda)$ satisfying $f(\lambda)/\lambda \to 0$ as $\lambda\to\infty$, see Appendix~\ref{appx:weyls-law}. The term $o(\lambda)$ collects contributions that grow strictly sublinearly in~$\lambda$. Hence, the slope of the smooth component of $N(\lambda)$ for large~$\lambda$ directly yields the area $ab$. In other words, the area appears as the \emph{slope} of the leading linear term in the asymptotic expansion of $N(\lambda)$. Since $N(\lambda)$ is computed directly from the ordered eigenvalues via $N(\lambda_k) = k$, the quantity $ab$ can be recovered by fitting a linear model to the tail of the sequence $\{(\lambda_k, k)\}_{k\gg 1}$. Indeed, rearranging~\eqref{eq:rectangle-weyl} gives
\[
\frac{N(\lambda)}{\lambda} \longrightarrow \frac{ab}{4\pi},
\qquad \lambda\to\infty,
\]
so that
\begin{equation}
\label{eq:area-limit}
ab = 4\pi \,\lim_{\lambda\to\infty} \frac{N(\lambda)}{\lambda}.
\end{equation}

Formula~\eqref{eq:area-limit} expresses the area purely in terms of the asymptotic density of the spectral sequence. No geometric information other than the eigenvalues themselves is required, and no mode labels $(m,n)$ enter the calculation. Thus the area of the rectangle is spectrally encoded in a particularly transparent way: it is the unique real number that normalizes the counting function~\eqref{eq:eigenvalue-counting-function} to unit slope in the high-frequency regime.

\section{Recovering the Side Lengths}
\label{sec:recovering-sides-fourier}
The individual side lengths $a$ and $b$ are encoded in the periodic fluctuations, that is, in the oscillatory remainder $F(\lambda)=N(\lambda)-N_{\text{smooth}}(\lambda)$ appearing in the decomposition~\eqref{eq:decomposition-smooth-osc}. As shown in Appendix~\ref{appx:periodic-orbit-length}, the Fourier transform of $F(\lambda)$ reveals sharp peaks at the periodic--orbit lengths of the underlying billiard. For a rectangular domain $\Omega_{a,b}=(0,a)\times(0,b)$ the shortest periodic classical trajectories are the bouncing--ball orbits parallel to the two sides, with lengths $2a$ and $2b$. Consequently, the length spectrum $S(L)$ defined in~\eqref{eq:spectral-density} exhibits dominant peaks at these two fundamental values. By locating these peaks, whose identification with the fundamental periodic–orbit lengths $2a$ and $2b$ is established in Appendix~\ref{appx:periodic-orbit-length}, we directly recover the side lengths:
\begin{equation}
\label{eq:side_length}
a=\frac{1}{2}L_{\text{short}},
\qquad
b=\frac{1}{2}L_{\text{long}},
\end{equation}
where $L_{\text{short}}$ and $L_{\text{long}}$ denote the two fundamental orbit lengths extracted from the Fourier spectrum. Besides the fundamental peaks at $2a$ and $2b$, the length spectrum typically contains additional peaks arising from longer periodic orbits,
\[
2ka,\qquad 2\ell b,\qquad 2\sqrt{(ka)^2 + (\ell b)^2},\qquad k,\ell\in\mathbb{N}.
\]
These additional peaks arise from the periodic--orbit length formula~\eqref{eq:periodic-orbit-length}. Depending on the spectral window, such harmonics may even exceed the amplitude of the fundamental ones. To resolve this ambiguity, we combine the Fourier data with the independently recovered area~$ab$: among all peak candidates $\{L_i\}$ we select the pair $(L_1,L_2)$ for which
\[
\frac{L_1}{2}\cdot\frac{L_2}{2}\approx ab.
\]
This geometric compatibility criterion reliably isolates the true lengths $2a$ and $2b$, even when their peaks are not dominant. Empirically, the recovered values converge to $2a$ and $2b$ as the number of available eigenvalues increases.

\section{Working Example}
To illustrate the reconstruction procedure, we provide a Mathematica notebook in the repository \cite{Sultanow2025}, see \textit{cone-operator-lab/mathematica/Hearing\_Rectangle\_Fourier.nb}. In this section we briefly walk through the key steps of the computation, using as a running example a rectangle with side lengths $a=1$ and $b=3$.

\medskip
To generate a sufficiently large sample of Dirichlet eigenvalues for the test rectangle with side lengths $a=1$ and $b=3$, we compute the values $\lambda_{m,n}=\pi^2(m^2/a^2+n^2/b^2)$ for $1\le m,n\le 800$ and extract the first $10\,000$ eigenvalues after sorting, see Listing~\ref{lst:1}.

\begin{lstlisting}[label={lst:1},language=Mathematica,caption={Generation of $\lambda_{m,n}$ on $\Omega_{1,3}$.}]
aTrue = 1.0;
bTrue = 3.0;
mMax = 800; nMax = 800;
ms = Range[1, mMax];
ns = Range[1, nMax];
lambdas = Table[Pi^2 (m^2/aTrue^2 + n^2/bTrue^2), {m, ms}, {n, ns}];
eigvals = Sort[Flatten[lambdas]][[6 ;; 10005]];
\end{lstlisting}

From the sorted list of eigenvalues we retain $10\,000$ values (after skipping the first five). For reference, the first twelve and the last three values in this list are

\[
\begin{array}{llll}
40.5750, & 43.8649, & 49.3480, & 49.3480, \\
57.0244, & 63.6041, & 66.8940, & 78.9568, \\
80.0535, & 89.9231, & 93.2129, & 93.2129, \\[4pt]
\multicolumn{4}{c}{\ldots} \\[4pt]
42430.5259, & 42440.3955, & 42440.3955, & 42440.3955
\end{array}
\]

These values serve as the input for the subsequent Weyl fit and Fourier–based reconstruction steps. To recover the area $ab$, we apply a linear Weyl fit to the upper asymptotic portion of the counting function, as implemented in Listing~\ref{lst:2}; this yields an estimate for the slope and thereby an estimate for the area.

\begin{lstlisting}[label={lst:2},language=Mathematica,caption={Weyl fit for estimating the area of $\Omega_{1,3}$.}]
startIndexA = Floor[0.3 Ntotal];
subVals = vals[[startIndexA ;;]];
subK = kVals[[startIndexA ;;]];
lm = LinearModelFit[Transpose[{subVals, subK}], x, x];
slopeA = lm["BestFitParameters"][[2]];
Ahat = 4 Pi slopeA;
\end{lstlisting}

To evaluate the length spectrum efficiently, we compute the quantity defined in~\eqref{eq:spectral-density} using a vectorized compiled function, as shown in Listing~\ref{lst:3}.

\begin{lstlisting}[label={lst:3},language=Mathematica,caption={Compiled evaluation of the length spectrum $S(L)$.}]
specFun =
  Compile[{{L, _Real}, {tKc, _Real, 1}, {FkWc, _Real, 1}},
   Module[{sRe = 0., sIm = 0., arg, c, s, lenLocal},
    lenLocal = Length[tKc];
    Do[
      arg = L*tKc[[k]];
      c = Cos[arg];
      s = Sin[arg];
      sRe += FkWc[[k]]*c;
      sIm += -FkWc[[k]]*s;
      ,
      {k, 1, lenLocal}
    ];
    sRe^2 + sIm^2
    ]
  ];
\end{lstlisting}

The dominant peaks of the length spectrum yield the fundamental periodic--orbit 
lengths $L_{\text{short}}\approx 1.9968$, $L_{\text{long}}\approx 6.0032$, corresponding to the theoretical values $2a=2$ and $2b=6$. Dividing by~$2$ gives the reconstructed side lengths $(\hat{a},\hat{b})\approx (0.9984,\;3.0016)$, which are accurate to within $0.2\%$. The derived aspect ratio is $\hat{R}=\hat{a}/\hat{b}\approx 0.3326$, in agreement with the true value $R = 1/3$.

\begin{figure}[H]
\includegraphics[width=0.9\textwidth]{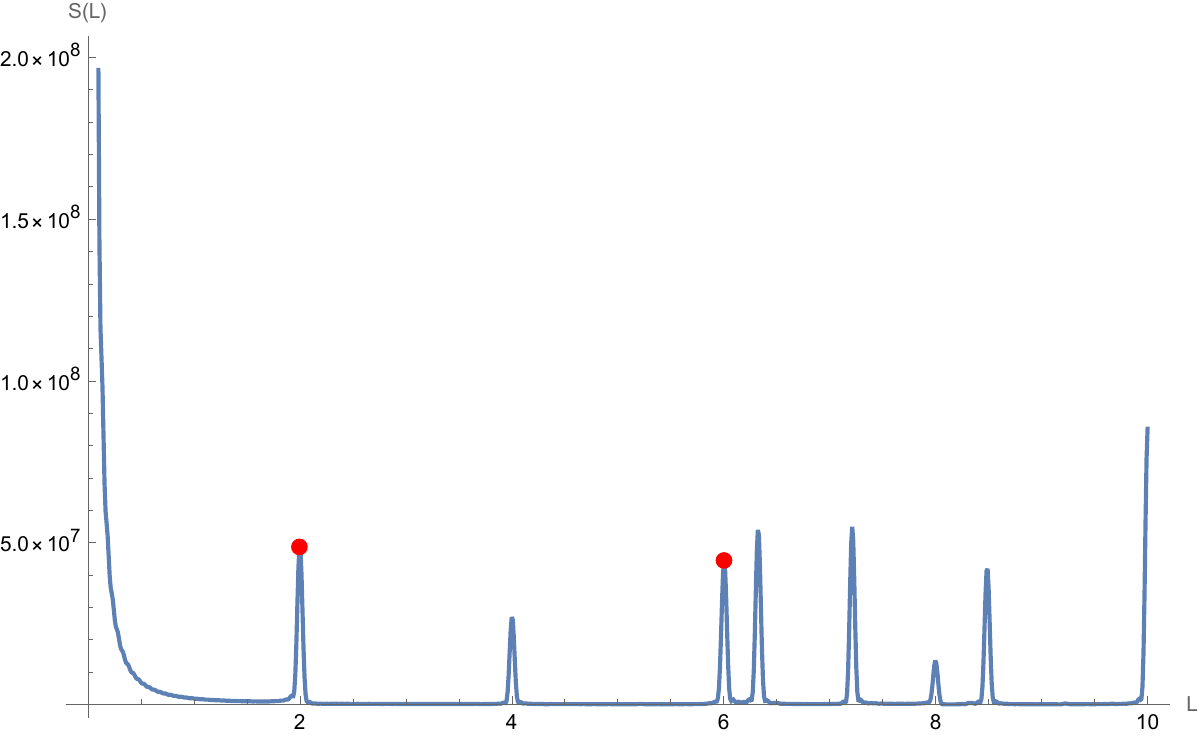}
\caption{Length spectrum $S(L)$ generated by equation~\eqref{eq:spectral-density}, with dominant peaks corresponding to the fundamental periodic--orbit lengths $2a$ and $2b$.}
\label{fig:peaks}
\end{figure}

\medskip
The peak plot in Figure~\ref{fig:peaks} shows these two dominant contributions clearly separated from the higher harmonics and composite orbit lengths. The two main peaks at $L\approx 2$ and $L\approx 6$ correspond exactly to the bouncing--ball orbits parallel to the sides of the rectangle, validating the periodic--orbit interpretation and confirming the accuracy of the reconstruction procedure.

\newpage
\section{Outlook and further Research}
The reconstruction of the side lengths of a rectangle from its Dirichlet spectrum is a first example of a broader inverse problem: detecting geometric symmetries directly from spectral data.\\

A natural intermediate step between planar rectangles and fully curved three-dimensional bodies is the Dirichlet spectrum of a three-dimensional rectangular box (parallelepiped) in $\mathbb{R}^3$. For the domain $\Omega_{a,b,c} = (0,a)\times(0,b)\times(0,c)$, the Dirichlet eigenvalues separate exactly as in the planar case, see Pinchover and Rubinstein~\cite[p.~357]{PinchoverRubinstein2005}:
\[
\lambda_{l,n,m}
= \Big(\tfrac{l\pi}{a}\Big)^2
+ \Big(\tfrac{n\pi}{b}\Big)^2
+ \Big(\tfrac{m\pi}{c}\Big)^2,
\qquad l,n,m\in\mathbb{N}
\]
The philosophy underlying the reconstruction of the side lengths of a rectangle extends naturally to the recovery of all three side lengths of a rectangular box. This provides a clean higher-dimensional analogue of our results and offers a stepping stone toward the spectral study of more complex three-dimensional geometries.

\bigskip
A particularly compelling next step is the spectral classification of convex bodies in $\mathbb{R}^3$ according to their symmetry type. Three fundamental classes are:

\begin{itemize}
  \item \textbf{Spheres (balls)}, determined by a single geometric parameter. Their Laplace spectra exhibit complete rotational symmetry and a highly structured multiplicity pattern.

  \item \textbf{Rotational ellipsoids}, determined by two distinct semi-axes. Their spectra remain separable in spheroidal coordinates, and the anisotropy is reflected in the splitting of formerly degenerate eigenvalue clusters.

  \item \textbf{General ellipsoids}, with three independent semi-axes. These exhibit no spectral degeneracy except for symmetry-imposed coincidences, and the fine structure of the ordered spectrum encodes the full anisotropy of the geometry.
\end{itemize}

Extending the presented techniques to three-dimensional bodies raises several mathematically rich questions: How stable is symmetry detection under spectral perturbations? To what extent can the semi-axes of an ellipsoid be recovered from a finite prefix of the spectrum? And can one characterize the exact spectral signatures that distinguish spherical, spheroidal, and fully triaxial geometries?

\section*{Acknowledgements}
I am grateful to Prof.\ Elmar Schrohe for suggesting the study of inverse spectral questions and for inspiring the investigation of how geometric information can be recovered from eigenvalue data.

\clearpage
\appendix

\section{Weyl’s Law}
\label{appx:weyls-law}
By the classical Weyl asymptotic law in Chavel's formulation \cite[p.~184]{Chavel2006} (see also Buser~\cite[p.~182]{Buser2010} and \cite[p.~18]{Levitin_etal2023}), the eigenvalue counting function~\eqref{eq:eigenvalue-counting-function} satisfies
\begin{equation}\label{eq:weyl-chavel}
N(\lambda)\;\sim\;
\frac{\omega_n\,\operatorname{Vol}(M)}{(2\pi)^n}\,
\lambda^{n/2},
\qquad \lambda\to\infty,
\end{equation}
where $\omega_n$ denotes the volume of the unit ball in $\mathbb{R}^n$. In dimension $n=2$, this becomes
\[
N(\lambda)\;\sim\;\frac{\omega_2\,\operatorname{Area}(M)}{(2\pi)^2}\,\lambda,
\qquad \lambda\to\infty.
\]
Since $\omega_2=\pi$ is the area of the unit disk, we obtain
\[
\frac{\omega_2}{(2\pi)^2}
= \frac{\pi}{4\pi^2}
= \frac{1}{4\pi},
\]
and therefore
\[
N(\lambda)
\sim \frac{\operatorname{Area}(M)}{4\pi}\,\lambda = \frac{ab}{4\pi}\,\lambda
\]
for a rectangle $M=\Omega_{a,b}$. As already noted by Levitin et al.\ \cite[p.~18]{Levitin_etal2023}, the Weyl term must be supplemented by a remainder term, denoted $R(\lambda)$, which in turn corresponds to de~Bruijn’s $o$–symbol. Its appearance follows from the standard equivalence $f(\lambda)\sim g(\lambda)\Longleftrightarrow f(\lambda)=g(\lambda)+o(g(\lambda))$. Applied to 
\[
g(\lambda)=\frac{\operatorname{Area}(M)}{4\pi}\,\lambda,
\]
this yields
\[
N(\lambda)
= \frac{\operatorname{Area}(M)}{4\pi}\,\lambda 
+ o\!\left(\frac{\operatorname{Area}(M)}{4\pi}\,\lambda\right)
= \frac{\operatorname{Area}(M)}{4\pi}\,\lambda + o(\lambda)
= \frac{ab}{4\pi}\,\lambda + o(\lambda).
\]
De~Bruijn’s $o$–symbol \cite[p.~3]{deBruijn1958} is understood in its standard asymptotic sense, denoting a negligible error term, i.e.\ a function $f(\lambda)$ satisfying $f(\lambda)=o(\lambda)$, equivalently $f(\lambda)/\lambda \to 0$ as $\lambda\to\infty$. De~Bruijn formulates the definition in terms of sequences $f(n)$, writing $f(n)/g(n)\to 0$ as $n\to\infty$; the continuous version used here is entirely analogous. Thus $o(\lambda)$ collects all contributions that grow strictly slower than $\lambda$ and are asymptotically negligible relative to the leading linear term.

\newpage
\section{Periodic--orbit length}
\label{appx:periodic-orbit-length}
In semiclassical periodic--orbit theory the oscillatory part of the spectral density is a superposition (i.e.\ a linear combination) of complex exponential terms of the form $\exp\!\big(i\,L_{\mathrm{po}}\sqrt{\lambda}\big)$, where $L_{\mathrm{po}}$ denotes the length of a periodic classical orbit. This structure follows from the fact that the classical action of an orbit satisfies $S_{M_1,M_2} = L_{M_1,M_2}\sqrt{E}$; see Brack and Bhaduri \cite[p.~95]{BrackBhaduri1997}. In their notation the semiclassical phase appears as $S_{M_1,M_2}/\hbar = (p/\hbar)\,L_{M_1,M_2}$, where $p=\sqrt{E}$ is the classical momentum. In our setting of Laplace eigenvalues $\lambda = E$ and dimensionless units with $\hbar = 1$, this reduces to the phase $L_{M_1,M_2}\sqrt{\lambda}$, which is precisely the exponent occurring in the Fourier factor $e^{-iL\sqrt{\lambda_k}}$. The real trace contributions arise from combining $\exp(\pm i L_{\mathrm{po}}\sqrt{\lambda})$ to give the customary cosine form. This motivates analyzing the discrete Fourier transform, which can also be found in the work of Marklof \cite[p.~177]{Marklof1998}:
\begin{equation}
\label{eq:spectral-density}
S(L) = \left| \sum_{k} F(\lambda_k)\, e^{-iL\sqrt{\lambda_k}} \right|^2,
\end{equation}
whose dominant peaks correspond to the fundamental periodic--orbit lengths of the billiard. For the rectangular billiard, Brack and Bhaduri \cite[p.~95]{BrackBhaduri1997} derive the explicit periodic--orbit length formula
\begin{equation}
\label{eq:periodic-orbit-length}
L_{M_1,M_2} = 2\sqrt{(M_1 a)^2 + (M_2 b)^2}.
\end{equation}
The fundamental orbits $(M_1,M_2) = (1,0)$ and $(0,1)$ therefore have lengths $2a$ and $2b$, respectively, and thus the shortest periodic trajectories in a rectangular billiard are of lengths $2a$ and $2b$. According to the geometric--spectral duality underlying Poisson summation and periodic--orbit theory (see, for example, the derivation of the periodic--orbit sum for the rectangular billiard in \cite[p.~96]{BrackBhaduri1997}), these fundamental orbit lengths (namely $2a$ and $2b$) appear as dominant peaks in the Fourier transform of the spectral fluctuations.

\clearpage
\bibliographystyle{unsrt}
\bibliography{main}

\end{document}